\documentclass[aps, amsmath, amssymb, aps, pre,twocolumn,groupedaddress,superscriptaddress,nofootinbib,longbibliography]{revtex4-1}

\usepackage{graphicx}
\usepackage[suffix=, outdir=./figs/]{epstopdf}
\usepackage[bottom]{footmisc}
\usepackage{enumitem}

\usepackage{array}
\usepackage{xfrac}
\usepackage[colorlinks,linkcolor=black,citecolor=blue,urlcolor=black]{hyperref}

\DeclareMathOperator{\LEs}{LE}

\newtheorem{theorem}{Theorem}
\newtheorem{lemma}{Lemma}
\newtheorem{corollary}{Corollary}

\begin{document}

\title{The Lyapunov dimension, convergency and entropy \\
for a dynamical model of Chua memristor circuit}

\author{G.A. Leonov}
\affiliation{Faculty of Mathematics and Mechanics, St. Petersburg State University,
Peterhof, St. Petersburg, Russia}
\author{N.V. Kuznetsov}
\email[]{Corresponding author: nkuznetsov239@gmail.com}
\affiliation{Faculty of Mathematics and Mechanics, St. Petersburg State University,
Peterhof, St. Petersburg, Russia}
\affiliation{Department of Mathematical Information Technology,
University of Jyv\"{a}skyl\"{a}, Jyv\"{a}skyl\"{a}, Finland}
\date{\today}


\begin{abstract}
\end{abstract}
\maketitle

\section*{Introduction}
For the study of chaotic dynamics and dimension of attractors
the concepts of the Lyapunov exponents \cite{Lyapunov-1892}
was found useful and became widely spread
\cite{GrassbergerP-1983,EckmannR-1985,ConstantinFT-1985,AbarbanelBST-1993}.
Such characteristics of chaotic behavior, as
the Lyapunov dimension \cite{KaplanY-1979}
and the entropy rate \cite{Kolmogorov-1959,Sinai-1959,AdlerKA-1965},
can be estimated via the Lyapunov exponents
\cite{Millionschikov-1976,Pesin-1977,KaplanY-1979,Young-2013}.
In this work an analytical approach to the study
of the Lyapunov dimension, convergency and entropy
for a dynamical model of Chua memristor circuit is demonstrated.

\section{A dynamical model of the Chua memristor circuit}

Consider one of the Chua memristor models \cite[eq.25]{CorintoF-2017}
\begin{equation}\label{chua-memristor}
\begin{aligned}
  & \dot x = \alpha(m_0-1)x+\alpha y-\alpha m_1x^3+\alpha x_0,  \\
  & \dot y = x-y+z, \ \dot z = \beta y-\gamma z
\end{aligned}
\end{equation}
with real parameters $\alpha,\beta,m_0,m_1,\gamma$, $x_0$,
and suppose that
\(
  \alpha m_1 > 0.
\)
For a survey on memristor circuits see, e.g. \cite{Tetzlaff-2014-book}.
System \eqref{chua-memristor} with $x_0=0$  describes the dynamics
of the Chua oscillator with cubic nonlinearity \cite{Chua-1993,HuangPWF-1996}.
If $\gamma=0$ and $x_0=0$, %$|x_0| >  \frac{2}{3}(m_0-1)\sqrt{\frac{m_0-1}{3m_1}}$
then system \eqref{sys:ode} has equilibria
$u^{eq}_0=(0,0,0)$ and
$u^{eq}_\pm = (\pm \sqrt{\frac{m0 - 1}{m_1}},0,\mp \sqrt{\frac{m0 - 1}{m_1}})$ for $m_0>1$.
Represent system \eqref{chua-memristor}
as an autonomous differential equation of general form:
\begin{equation}\label{sys:ode}
  \dot{u} = f({u}),
\end{equation}
where $u=(x,y,z) \in U=\mathbb{R}^3$ and the
continuously differentiable vector-function $f: \mathbb{R}^3 \to \mathbb{R}^3$
is the right-hand side of system \eqref{chua-memristor}.
Define by $u(t,u_0)$ a solution of \eqref{sys:ode} such that $u(0,u_0)=u_0$,
and consider the evolutionary operator $\varphi^t(u_0) = u(t,u_0)$.
We assume the uniqueness and existence of solutions of \eqref{sys:ode} for $t \in [0,+\infty)$.
Then system \eqref{sys:ode} generates a dynamical system
$\{\varphi^t\}_{t\geq0}$.
Let a nonempty set $K \subset U$ be invariant with respect to $\{\varphi^t\}_{t\geq0}$,
i.e. $\varphi^t(K) = K$ for all $t \geq 0$.
For example, as the set $K$ one can consider various types of attractors
of system \eqref{sys:ode} (see, e.g. examples from \cite{HuangPWF-1996,BilottaP-2008}).
Recently the \emph{classification of local attractors as being hidden or self-excited}
was introduced in connection with the discovery of the first hidden attractor
in the classical Chua model with the saturation nonlinearity
\cite{KuznetsovLV-2010-IFAC,LeonovKV-2011-PLA,BraginVKL-2011,LeonovKV-2012-PhysD,KuznetsovKLV-2013,KiselevaKKKLYY-2017,StankevichKLC-2017}:
an attractor is called a \emph{self-excited attractor}
if its basin of attraction
intersects with any open neighborhood of an equilibrium,
otherwise, it is called a \emph{hidden attractor} \cite{LeonovKV-2011-PLA,LeonovK-2013-IJBC,LeonovKM-2015-EPJST,Kuznetsov-2016}. %,
For example, hidden attractors can be found in various memristive circuits,
see, e.g. \cite{PhamVVJKH-2016,HuLLYZ-2017-HA,BaoBWCX-2017-HA,BaoBLWW-2016-HA,PhamVVHV-2016-HA,Vaidyanathan-2017-SCI-HA,VaidyanathanPV-2017-HA,PhamVVTT-2017-HA,PhamVVWH-2017-HA,RochaRK-2017-HA,SahaSRC-2015-HA,SemenovKASVA-2015,PhamVJWV-2014-HA,PhamVVLV-2015-HA,ChenYB-2015-HA-EL,ChenLYBXW-2015-HA}.
Remark that hidden attractors are not connected with equilibria and, thus,
do not related with the Shilnikov scenario of chaos \cite{AfraimovichGLShT-2014}.
In the works \citep{BianchiKLYY-2015,BlagovKLYY-2015,LeonovKYY-2015-TCAS,KuznetsovLYY-2017-CNSNS}
it is demonstrated the difficulties of reliable simulation
of the phase-locked loops circuits in SPICE and MATLAB Simulink,
caused by hidden attractors with narrow basins of attraction.

%For $x_0=0$ system \eqref{sys:lorenz-general}
%has the following equilibria: $u^{eq}_0 = (0,0,0)$ and $u^{eq}_{\pm}=$
%\[
%  \bigg
%  (\pm \tfrac{\sqrt{\beta m_0 - \beta + \gamma m_0}}{\sqrt{m_1(\beta + \gamma)}},
%  \pm \tfrac{\gamma\sqrt{\beta m_0 - \beta + \gamma m_0}}{\sqrt{m_1(\beta + \gamma)^3}},
%  \mp \tfrac{\beta\sqrt{\beta m_0 - \beta + \gamma m_0}}{\sqrt{m_1(\beta + \gamma)^3}}
%  \bigg)
%\]
%for $\beta m_0 - \beta + \gamma m_0>0$ and $\beta\sqrt{\beta m_0 - \beta + \gamma m_0>0$.

Consider linearization of system \eqref{sys:ode}
along the solution $u(t,u_0)=\varphi^t(u)$:
\begin{equation} \label{sfl}
  \begin{aligned}
    & \dot v = J(\varphi^t(u))v,
    \quad J(u) = Df(u),
  \end{aligned}
\end{equation}
where $J(u)$ is the $3\!\times\!3$ Jacobian matrix
\[
  J(u) =\left(
       \begin{array}{ccc}
         \alpha(m_0-1)-3\alpha m_1x^2 & \alpha & 0 \\
         1 & -1 & 1 \\
         0 & \beta & -\gamma \\
       \end{array}
     \right)
\]
and it can be represented as
\(
  J=J(0)-3\alpha m_1x^2I_1
\) with
\[
J_0= J(0) =\left(
       \begin{array}{ccc}
         \alpha(m_0-1)&\alpha & 0 \\
         1 & -1 & 1 \\
         0 & \beta & -\gamma \\
       \end{array}
     \right), \
I_1 =\left(
       \begin{array}{ccc}
         1&0 & 0 \\
         0 & 0 & 0 \\
         0 & 0 & 0 \\
       \end{array}
     \right).
\]

Let for any $t > 0$ and any $u \in U$ the ordered sequence
$\lambda_1(u) \ge\!\cdots\!\ge \lambda_n(u)$, where
$\lambda_i(u) = \lambda_i(\tfrac{1}{2} (J(u) + J(u)^{*})$,
$i = 1,\dots,n$ be the eigenvalues of the symmetrized Jacobian matrix
\(
%\begin{equation} \label{SJS}
  \tfrac{1}{2} \left(J(u) + J(u)^{*} \right).
%\end{equation}
\)
\begin{lemma}
\(
 \lambda_j(0) \geq \lambda_j(u), \quad j=1,2,3
\)
\end{lemma}

Then from Corollary~\ref{thm:dLKYeig} (see in the Appendix) we get

\begin{theorem}
\(
 d_{\rm L}^{\rm KY}\big(\{\lambda_{j}(u)\}_{i=1}^3\big)
 \leq
 d_{\rm L}^{\rm KY}\big(\{\lambda_{j}(0)\}_{i=1}^3\big)
\)
\end{theorem}

If $J(0)$ have simple real eigenvalues $\lambda_i(J(0))$, then
$\lambda_i(J(0)) = \lambda_i(\tfrac{1}{2}(J(0)+J^*(0)))$ and
we get the following result

\begin{corollary}
Let $u^{eq}_0=(0,0,0)$ be one of the equilibria of system \eqref{chua-memristor}
and the matrix $J(0)$ have simple real eigenvalues.
Then the exact Lyapunov dimension
of any compact invariant set $K \ni u^{eq}_0$ is defined as
\[
 \dim_{\rm L}K =
 d_{\rm L}^{\rm KY}\big(\{\lambda_{j}(0)\}_{i=1}^3\big).
\]
\end{corollary}

By Theorem~\ref{theorem:l1l2} (see in the Appendix) we get
\begin{theorem}
If
\(
 \lambda_{1}(0)+\lambda_{2}(0) <0,
\)
then any bounded solution of system \eqref{sys:ode}
tends to the stationary set.
\end{theorem}

%For example, consider \cite{CorintoF-2017}
%system \eqref{chua-memristor} with parameters
%$\beta = 16$, $m_0 = 1-0.143$, $m_1 = 1$, $\alpha=10$, and $x_0=0$.
%Then \eqref{chua-memristor} has equilibria
%$u^{eq}_0=(0,0,0)$ and
%$u^{eq}_\pm = (\pm 3/2,0, \mp 3/2)$.
%
%Let $m_1=1$, $c=1-m_0$, $\alpha=a$ and $\beta=-b$
%\[
%\begin{aligned}
%  & \dot x = -acx+ay-ax^3,  \\
%  & \dot y = x-y+z, \\
%  & \dot x = -by
%\end{aligned}
%\]
%Characteristic polynomial is
%\[
%  t = {\rm Tr}(J(u)) = - 3ax^2 - ac - 1
%\]
%\[
%  \det(pI-J(u)) =
%    p^3 + (3ax^2 + ac + 1)p^2 + (3ax^2 - a + b + ac)p + 3abx^2 + abc
%\]
%\[
%  \det(pI-J(u)) =p^3 -t p^2 + (-t-1-a+b)p + b(-t-1)
%\]
%\[
%p_1(x)+p_2(x)+p_3(x) = t = - 3ax^2 - ac - 1
%\]

\newpage
\section*{Appendix. Exact and finite-time Lyapunov dimension}
Suppose that $\det J(u) \neq 0 \quad \forall u \in U$.
Consider a fundamental matrix of linearized system \eqref{sfl}
$D\varphi^t(u)$ such that $D\varphi^0(u) = I$, where $I$ is a unit $3 \times 3$ matrix.
Let $\sigma_i(t,u) = \sigma_i(D\varphi^t(u))$, $i = 1,2,3$,
be the singular values of $D\varphi^t(u)$
with respect to their algebraic multiplicity
ordered so that $\sigma_1(t,u) \geq \sigma_2(t,u) \geq \sigma_3(t,u) > 0$
for any $u \in U$ and $t \geq 0$.
Consider the set of \emph{the finite-time Lyapunov exponents} at the point $u_0$ $\{\LEs_i(t,u_0)=\frac{1}{t}\ln\sigma_{j}(t,\!u_0)\}_{i=1}^3$
ordered by decreasing  for $t > 0$.

Introduce the \emph{Kaplan-Yorke formula \cite{KaplanY-1979} with respect
to the ordered set} $\lambda_1\geq... \geq \lambda_n$:
\begin{equation}\label{lftKY}
   d^{\rm KY}(\{\lambda_i\}_{i=1}^3)\!=\!
   j+\frac{\sum_{i=1}^{j}\lambda_i}{|\lambda{j + 1}|},  \
   j\!=\!\max\{m\!: \sum_{i=1}^{m}\!\lambda_i\!\geq\!0\},
\end{equation}
where $d^{\rm KY}(\{\lambda_i\}_{i=1}^3)=0$ for $j=0$ and
$d^{\rm KY}(\{\lambda_i\}_{i=1}^3)=3$ for $j=3$.
Then
the \emph{finite-time local Lyapunov dimension}
\cite{Kuznetsov-2016-PLA,KuznetsovLMPS-2017}
at a certain point $u_0$ can be defined as
\[
   \dim_{\rm L}(t,u_0) = d^{\rm KY}(\{\LEs_i(t,u_0)\}_{i=1}^3).
\]
and the \emph{finite-time Lyapunov dimension}
of invariant closed bounded set $K$
is as follows
\begin{equation}\label{DOmaptmax}
  \dim_{\rm L}(t, K) = \sup\limits_{u_0 \in K} \dim_{\rm L}(t,u_0).
\end{equation}
In this approach the use of Kaplan-Yorke formula \eqref{lftKY} with
the finite-time Lyapunov exponents is justified
by the \emph{Douady--Oesterl\'{e} theorem} \cite{DouadyO-1980},
which implies that for any fixed $t > 0$
the Lyapunov dimension of the map $\varphi^t$ with respect
to a closed bounded invariant set $K$, defined by \eqref{DOmaptmax},
is an upper estimate of the Hausdorff dimension of the set $K$:
\(
  \dim_{\rm H}K \leq \dim_{\rm L}(t, K).
\)
For the estimation of the Hausdorff dimension of invariant closed bounded set $K$
one can use the map $\varphi^t$ with any time $t$
(e.g. $t=0$ leads to the trivial estimate $\dim_{\rm H}K \leq 3$)
and, thus, the best estimation is
\(
  \dim_{\rm H}{K} \le \inf_{t\geq0}\dim_{\rm L} (t, K).
\)
The following property
\begin{equation}\label{DOlim}
  \inf_{t\geq0}\sup\limits_{u \in K} \dim_{\rm L}(t,u)
  = \liminf_{t \to +\infty}\sup\limits_{u \in K} \dim_{\rm L}(t,u)
\end{equation}
allows one to introduce the \emph{Lyapunov dimension}  %of $K$ as
\cite{Kuznetsov-2016-PLA}
\begin{equation}\label{DOinf}
  \dim_{\rm L} K
  = \liminf_{t \to +\infty}\sup\limits_{u \in K} \dim_{\rm L}(t,u).
\end{equation}
If the maximum of local Lyapunov dimensions on the global attractors,
which involves all equilibria, is achieved at an equilibrium point $u^{cr}_{eq}$,
i.e. $\dim_{\rm L} u^{cr}_{eq} = \max_{u_0 \in K} \dim_{\rm L} u_0$,
then this allows one to get the \emph{exact Lyapunov dimension}.
(this term was suggested by Doering~et~al.~in~\cite{DoeringGHN-1987}).
In general, a \emph{conjecture on the Lyapunov dimension of self-excited attractor} \cite{KuznetsovLMPS-2017}
is that for a typical system
the Lyapunov dimension of a self-excited attractor
does not exceed the Lyapunov dimension of one of unstable equilibria,
the unstable manifold of which intersects with the basin of attraction
and visualize the attractor.

%Recall that a set $K$ with noninteger Hausdorff dimension
%is referred to as a \emph{fractal set} \cite{EckmannR-1985}
%and, when such set $K$ is an attractor,
%it is called a \emph{strange attractor}
%\cite{RuelleT-1971,GrebogiOY-1987}.

In contrast to the finite-time Lyapunov dimension \eqref{DOmaptmax},
the Lyapunov dimension \eqref{DOinf}
is \emph{invariant under smooth change of coordinates}
\cite{KuznetsovAL-2016,Kuznetsov-2016-PLA}.
This property and a proper choice of smooth change of coordinates
may significantly simplify the estimation of the Lyapunov dimension
of dynamical system.

%Based on this property and using a proper choice of smooth change of coordinates,
%in many cases it is possible to compute the exact Lyapunov dimension analytically
%by the method,  prosed by Leonov \cite{Leonov-1991-Vest,LeonovB-1992,Leonov-2012-PMM,Kuznetsov-2016-PLA}.

Consider an effective analytical approach, proposed by Leonov
\cite{Leonov-1991-Vest,LeonovB-1992,Kuznetsov-2016-PLA}.
Let for any $t > 0$ and any $u_0 \in U$ the ordered sequence
$\lambda_1(u_0, S) \ge\!\cdots\!\ge \lambda_n(u_0, S)$, where
$\lambda_i(u_0, S) = \lambda_i(S \varphi^t(u_0))$, $i = 1,\dots,n$
be the eigenvalues of the symmetrized Jacobian matrix
\begin{equation} \label{SJS}
  \frac{1}{2} \left( S J(\varphi^t(u_0)) S^{-1} + (S J(\varphi^t(u_0)) S^{-1})^{*}\right).
\end{equation}

\begin{theorem}\label{theorem:th1}
If there exist an integer $j \in \{1,\ldots,n-1\}$,
a real $s \in [0,1]$, a differentiable scalar function $V: U \subseteq \mathbb{R}^n \to \mathbb{R}^1$, and a nonsingular $n\times n$ matrix $S$
such that
\begin{equation}\label{ineq:weilSVct}
  \sup_{u \in U} \big( \lambda_1 (u,S) + \cdots + \lambda_j (u,S)
  + s\lambda_{j+1}(u,S) + \dot{V}(u) \big) < 0,
\end{equation}
where $\dot{V} (u) = ({\rm grad}(V))^{*}f(u)$,
then for a compact invariant set $K\subset U$ we have
\[
   \dim_{\rm H}K \leq
    \dim_{\rm L}(\{\varphi^t\}_{t\geq0},K)
   < j+s.
\]
\end{theorem}
%This method is based on a combination of the Douady-Oesterl\'{e} approach
%with the direct Lyapunov method and
In the work \cite{Kuznetsov-2016-PLA}
it is shown how the method can be justified
by the invariance of the Lyapunov dimension of compact invariant set
with respect to the special smooth change of variables $h$
with $Dh(u)=e^{V(u)(j+s)^{-1}}S$, where $V$
is a differentiable scalar function and $S$ is a nonsingular $n \times n$ matrix.
For $S=0$ and $V(u) \equiv 0$ we have
\begin{corollary}[\cite{DouadyO-1980,Smith-1986,Leonov-1991-Vest,Kuznetsov-2016-PLA}]
\label{thm:dLKYeig}
%For a compact invariant set $K$ %we have
%\begin{equation}\label{dLKYeig}
\[
 \dim_{\rm H}K\!\leq\!\dim_{\rm L}K \!\leq\!
 \sup_{u \in K}d_{\rm L}^{\rm KY}\!\big(\{\lambda_{j}(u)\}_{i=1}^3\big).
\]
\end{corollary}

The following result~\cite{Leonov-2012-PMM} is useful for the study of global convergency.
\begin{theorem}\label{theorem:l1l2}
If there exist a continuously differentiable scalar function
$V: U \subseteq \mathbb{R}^n \to \mathbb{R}^1$ and
a non-degenerate $n\times n$ matrix $S$ exist such that
\begin{equation}\label{pmm36}
  \sup_{u \in U} \big(
  \lambda_1(u,S)+\lambda_2(u,S)+\dot V(u)
  \big)
  <0,
\end{equation}
then any bounded solution of system \eqref{sys:ode}
with any initial data $u_0 \in U$
tends to the stationary set of dynamical system $\{\varphi^t\}_{t\geq0}$ as $t\to+\infty$.
\end{theorem}
Remark that the stationary set can have any structure,
e.g. %can consists of a finite number of stable and unstable point or
be a line of  equilibria.

In \cite{BoichenkoL-1998,PogromskyM-2011} it is demonstrated how the above technique
can be effectively used to a derive constructive upper bound of
the sum of positive Lyapunov exponents and
the topological entropy \cite{AdlerKA-1965}
(the topological entropy is an analogue of the entropy defined earlier
by Kolmogorov and Sinai \cite{Kolmogorov-1959,Sinai-1959}).

\section*{Acknowledgment}
\vspace{-0.5cm}
We would like to thank Leon Chua, Fernando Corinto, and
Ronald Tetzlaff to draw our attention to the memristors
and for fruitful discussions.

\vspace{-0.5cm}
%\bibliographystyle{elsarticle-num}
%\bibliography{C:/Dropbox/bib/bib_nk,C:/Dropbox/bib/bib_leonov,C:/Dropbox/bib/bib_full,C:/Dropbox/bib/genlorenz-bib,C:/Dropbox/bib/bib-hidden}

\end{document}